\documentclass[oneside,english]{amsart}
\usepackage[T1]{fontenc}
\usepackage[latin9]{inputenc}
\usepackage{amsthm}
\usepackage{amssymb}
\usepackage{graphicx}
\usepackage{epsfig}
\usepackage{epstopdf}
\usepackage{graphics}
\usepackage{graphicx}

\makeatletter
\numberwithin{equation}{section}
\numberwithin{figure}{section}
\theoremstyle{plain}
\newtheorem{thm}{Theorem}
  \theoremstyle{plain}
  \newtheorem{lem}[thm]{Lemma}
  \newtheorem{prop}[thm]{Proposition}
    \newtheorem{corollary}[thm]{Corollary}
  \theoremstyle{plain}
  \newtheorem{rem}[thm]{Remark}

\usepackage{latexsym}

\makeatletter
\AtBeginDocument{
  
}

\makeatother

\AtBeginDocument{

}

\AtBeginDocument{
  
}

\makeatother

\usepackage{babel}

\begin{document}

\title[Exponential growth rates of Baumslag-Solitar groups]
{Minimal exponential growth rates of metabelian Baumslag-Solitar groups and lamplighter groups}
\begin{abstract}
We prove that for any prime $p\geq 3$ the minimal exponential growth rate of the Baumslag-Solitar group $BS(1,p)$ and the lamplighter group $\mathcal{L}_p=(\mathbb{Z}/p\mathbb{Z})\wr \mathbb{Z}$ are equal.
We also show that for $p=2$ this claim is not true and the growth rate of $BS(1,2)$ is equal to the positive root of $x^3-x^2-2$, whilst the one of the lamplighter group $\mathcal{L}_2$ is equal to the golden ratio $(1+\sqrt5)/2$.
The latter value also serves to show that the lower bound of A.Mann from \cite{Mann} for the growth rates of non-semidirect HNN extensions is optimal.
\end{abstract}
\author{Michelle Bucher, Alexey Talambutsa}
\thanks{Michelle Bucher is supported by Swiss National Science Foundation
project PP00P2-128309/1. Alexey Talambutsa is supported by Russian Science Foundation, project 14-50-00005.}
\maketitle

\section{Introduction}
Let $G$ be a finitely generated group. For any finite generating set $S$ of $G$ we can consider the \emph{exponential growth rate} of $G$ with respect to $S$ which is defined as follows: Any element $g\in G$ can be written as a finite product of elements in $S\cup S^{-1}$ and we define the length  $\ell_{G,S}(g)$ of $g$ as the minimum length of such a product. The growth function $F_{G,S}(n)$ counts the number of elements in a ball of radius $n$ centered at the identity, that is the number of elements $g\in G$ for which $\ell_{G,S}(g)\leqslant n$. Finally the \emph{exponential growth rate} of $G$ with respect to $S$ is the limit
$$\omega(G,S)=\lim\limits_{n\to \infty}(F_{G,S}(n))^{\frac1n}\geq 1.$$
Note that this limit always exists by submultiplicativity of the growth function (see \cite[VI.C.56]{dlHBook}).

The exponential growth rate $\omega(G,S)$ clearly depends on the choice of the generating set $S$ and one obtains a group invariant by considering the infimum over all finite generating sets:
\begin{equation}
\Omega(G)= \inf_{|S|<\infty} \{\omega(G,S)\}.
\label{growth-infimum}
\end{equation}

It is now natural to ask if there exists generating sets $S$ for which the equality $\Omega(G)=\omega(G,S)$ is realized. For the free group $\mathbb{F}_n$ of rank $n$, Gromov remarked in \cite[Example 5.13]{Gromov} that $\Omega(\mathbb{F}_n)$ is exactly $2n-1$ and is realized on any free generating set (with $n$ elements). Except for this example, very few exact values for $\Omega(G)$ have been computed. Known cases include free products $\mathbb{Z}_2*\mathbb{Z}_{p^k}$ \cite{Talambutsa2011} (the cases $p^k=3,4$ were proven earlier in \cite{Mann}), the free product $\mathbb{Z}_2*(\mathbb{Z}_2\times \mathbb{Z}_2)$ and the Coxeter group $\mathrm{PGL}(2,\mathbb{Z})$ \cite{Bucher-Talambutsa-IJM} and a few more examples in the references \cite{Bucher-Talambutsa-IJM,Talambutsa2011,Mann}. But the question of de la Harpe and Grigorchuk whether $\Omega(\pi_1(\Sigma_g))$ is realized on the canonical generators of the fundamental group of a closed surface $\Sigma_g$ with $g\geq 2$ is still open (see \cite[p.55]{Grigorchuk-dlHarpe}). While in many cases, the value $\omega(G,S)$ can be computed for some particular generating set $S$, it is usually much harder to find a generating set $S$ such that $\Omega(G)=\omega(G,S)$ and sometimes even impossible due to the existence of groups for which the infimum in \eqref{growth-infimum} is not attained (see \cite{Sambusetti, Wilson}). 

We consider two classes of metabelian groups: Baumslag-Solitar groups $BS(1,n)$ and lamplighter groups $\mathcal{L}_n=(\mathbb{Z}/n\mathbb{Z})\wr \mathbb{Z}$.
The growth functions of the Baumslag-Solitar groups
\begin{equation}
BS(1,n)=\langle a,t\mid tat^{-1}=a^n \rangle
\label{Baumslag-Solitar-presentation}
\end{equation}
with respect to the canonical generating set $S=\{a,t\}$ were computed by Collins, Edjvet and Gill in \cite{Collins-Edjvet-Gill}. The restricted wreath products $\mathcal{L}_n=(\mathbb{Z}/n\mathbb{Z})\wr \mathbb{Z}$ can be presented as
\begin{equation}
\mathcal{L}_n=\langle a,t \mid a^n=1,\,\,\, [t^kat^{-k},a]=1 \,\, (k=1,2,\ldots) \rangle.
\label{Baumslag-Solitar-presentation}
\end{equation}
To compute the growth function of $\mathcal{L}_n$ with respect to the set $\{a,t\}$
one can use formulas given by Parry in \cite{Parry}. Even though the formulas for the growth functions of $BS(1,n)$ and $\mathcal{L}_n$  were obtained by completely different methods and by use of different properties of the groups, we find that remarkably for all odd $n=2k+1$
\begin{equation}
\omega(BS(1,n),\{a,t\})=\omega(\mathcal{L}_n,\{a,t\})=\omega_k,
\label{GrowthEquality}
\end{equation}
where $\omega_k$ is the unique positive root of
\[T_k(x)=x^{k+1}-x^k-2x^{k-1}-\ldots-2x-2,\]
for $k\geq 1$. This is easily deduced from \cite{Parry} and \cite{Collins-Edjvet-Gill} in Lemma \ref{LemmaGrowthRatesOfBS}. Interestingly, this equality never holds for even $n$. We will see the case $n=2$ in more details.

Some inference for the equality \eqref{GrowthEquality} can be seen in the actions of the groups $BS(1,n)$ and $\mathcal{L}_n$ on their corresponding Bass-Serre trees. There is indeed a very strong similarity between these actions, which we exploit to prove the main result of the paper:

\begin{thm} Let $p$ be a prime. The minimal growth rate of the Baumslag-Solitar group $BS(1,p)$ and lamplighter groups $\mathcal{L}_p$ are realized on the canonical generators $\{a,t\}$:
\[
\Omega(\mathcal{L}_p)=\Omega(BS(1,p))= \omega_k ,\quad \mathrm{for \ } p=2k+1,
\]
$$\Omega(\mathcal{L}_2)=\frac{1+\sqrt5}2 <\Omega(BS(1,2))= \beta,$$
where $\beta\sim 1.69572$ is the unique positive root of $z^3-z^2-2$.
\label{Main Theorem}
\end{thm}

The exact computation $\Omega(\mathcal{L}_2)=(1+\sqrt5)/2$ gives a positive answer to the question of Mann \cite{Mann} whether the lower bound $\Omega(G)\geq (1+\sqrt{5})/2$ can be realized on a non-semidirect HNN extension. (The fact that  $\mathcal{L}_{2}$ is indeed a non-semidirect HNN extension will be shown in Section \ref{section: lamplighter groups}). Note that it follows from Theorem~\ref{Main Theorem} that this lower bound could never be realized on any of the Baumslag-Solitar groups $\Omega(BS(1,n))$ also for arbitrary integers $n\geq 2$.

The lower bounds for the growth rates in Theorem \ref{Main Theorem} are obtained by looking at the actions on the corresponding Bass-Serre tree, finding free submonoids using a local variant of the classical ping-pong lemma (Lemma \ref{Key Lemma} here) and computing their growth with Lemma \ref{lem: lengthMonoid}. Interestingly, all the minimal growth rates are in fact realized as the growth rate of some free submonoid. The Bass-Serre trees of $ \mathcal{L}_p$ and $BS(1,p)$ are both $(p+1)$-regular trees, but the corresponding actions are of course different. Nevertheless, when $p$ is odd, the same method applies to give the lower bound of Theorem \ref{Main Theorem}, which we abstract in the following theorem:

\begin{thm}\label{Main Theorem 2} Let $G=H*_\theta$ be an HNN extension relative to an isomorphism $\theta:A\rightarrow B$ with $A=H$ and $B$ a normal subgroup of prime index $p$ in $H$. Then
$$\Omega(G)\geq \frac{1+\sqrt{5}}{2}, \quad \mathrm{for \ } p=2,$$
$$\Omega(G) \geq \omega_k, \quad \mathrm{for \ } p=2k+1.$$
\end{thm}

Together with the equalities \eqref{GrowthEquality} proven in Lemma  \ref{LemmaGrowthRatesOfBS} this immediately implies Theorem \ref{Main Theorem}, except in the case of $BS(1,2)$. For this last group, a finer analysis of its action on its Bass-Serre tree will be needed.

The question of Mann mentioned above was prompted by his proof of the lower bound $\Omega(G)\geq (1+\sqrt{5})/2$ for any non-semidirect HNN extension $G$ (see \cite{Mann}), using the cute algebraic observation that a hyperbolic element and a nontrivial conjugate of it generate a free monoid with growth rate equal to the golden ratio. Our proof for the case $p=2$ of Theorem \ref{Main Theorem 2} also holds for any non-semidirect HNN extension and gives an alternative geometric proof to Mann's inequality.

Finally, as an application of Theorem \ref{Main Theorem}, we can compute the minimal growth rate of the wreath product $\mathbb{Z}\wr \mathbb{Z}$. Indeed,  as was already noted by Shukhov  in \cite{Shukhov},
one can deduce from \cite{Collins-Edjvet-Gill} that
\begin{equation}
\lim_{n\to \infty}\omega({BS(1,n),\{a,t\}})=1+\sqrt2.
\label{LimitValue}
\end{equation}
Since the wreath product $\mathbb Z \wr \mathbb Z$ can be viewed as an extension of the groups $\mathcal{L}_p$, combining Theorem \ref{Main Theorem} and Parry's computations for $\mathbb Z \wr \mathbb Z$, we obtain

\begin{corollary} The minimal growth rate of the restricted wreath product $\mathbb Z \wr \mathbb Z=\langle a,t \mid [t^kat^{-k},a]=1 \,\, (k=1,2,\ldots) \rangle$ is realized on the set $\{a,t\}$ and
$$\Omega(\mathbb Z \wr \mathbb Z)=\omega(\mathbb Z \wr \mathbb Z,\{a,t\})=1+\sqrt{2}.$$
\label{cor: ZwrZ}
\end{corollary}

\subsection*{Acknowledgements} We thank Murray Elder for helpful discussions in the preparation of this work and Tatiana Smirnova-Nagnibeda for pointing out some useful references.
\section{Bass-Serre tree for an HNN extension}\label{section HNN}

Let $G=H*_\theta$ be the HNN extension of $H$ relative to the isomorphism $\theta: A\rightarrow B$ between the two subgroups $A,B$ of $H$. Following \cite{Mann} we call $H*_\theta$ a \emph{non-semidirect} HNN-extension if at least one of the subgroups $A$ or $B$ is a proper subgroup in $H$. If $H=\langle S_H\mid R_H \rangle $ is a presentation of $H$, then $G$ admits the presentation
$$G=\langle S_H, t \mid R_H, tat^{-1}=\theta(a) \ \forall a\in A\rangle.$$
There is a natural surjection $\varphi:G\rightarrow \mathbb{Z}$ defined by sending the generators $S_H$ to $0$ and $t$ to $1$.

The vertices of the associated Bass-Serre tree $T$ of $G$ are the right cosets of $G$ by $H$ and the edges are the right cosets of $G$ by $B$,
$$T^0=G/H, \quad T^1=G/B.$$
The edge $gB\in T^1$ has vertices $gH$ and $gtH$. This is a tree of valency $[H:A]+[H:B]$. The group $G$ acts on $T$ by left multiplication.

Since the natural surjection  $\varphi:G\rightarrow \mathbb{Z}$ is trivial on $H$, it induces a map $\overline{\varphi}:T^0\rightarrow \mathbb{Z}$ which sends vertices $v,w$ of an edge of $T^1$ to images satisfying $|\overline{\varphi}(v)-\overline{\varphi}(w)|=1$. This allows us to define an orientation on the edges by giving an edge from $v$ to $w$ with $\overline{\varphi}(w)-\overline{\varphi}(v)=1$ the positive orientation. This allows us to distinguish between two types of neighboors to any vertex $v$: the $[H:A]$ vertices $w$ such that $\overline{\varphi}(w)=\overline{\varphi}(v)-1$ which we call the {\it direct ascendants} of $v$, and the $[H:B]$ vertices $w$ such that $\overline{\varphi}(w)=\overline{\varphi}(v)+1$, which we call the {\it direct descendants} of $v$. We further call a vertex $z$ a {\it ascendant}, respectively an {\it descendant}, of $v$ if there is a sequence $v=w_0,w_1,\dots,w_\ell=z$ such that $w_i$ is a direct ascendant, resp. direct descendant, of $w_{i-1}$ for $1\leq i\leq \ell$. In our examples, $[H:A]=1$, which means that there is only one direct ascendant to any vertex. We will also use the terminology that a vertex $v$ is {\it above}, respectively {\it below}, a vertex $w$ if $v$ is an ascendant, resp. descendant, of $w$.

Since the action of $G$ on $T$ preserves the orientation on the edges defined above, it is immediate  that $G$ acts on $T$ without inversions. Thus there are two types of elements: elliptic elements $g\in G$ have a fixed point on $T$ and are thus conjugated to $H$, and hyperbolic elements $g\in G$ have no fixed point and possess a unique invariant geodesic $L_g$, called the axis of $g$, on which $g$ acts by translation. Note that any element $g\in G$ which is not in the kernel of $\varphi:G\rightarrow \mathbb{Z}$ necessarily is hyperbolic, so in particular, any generating set of $G$ contains a hyperbolic element. Such hyperbolic elements will be called positive, respectively negative according to their image in $\mathbb{Z}$ being positive or negative.

\begin{figure}
\includegraphics[width=80mm]{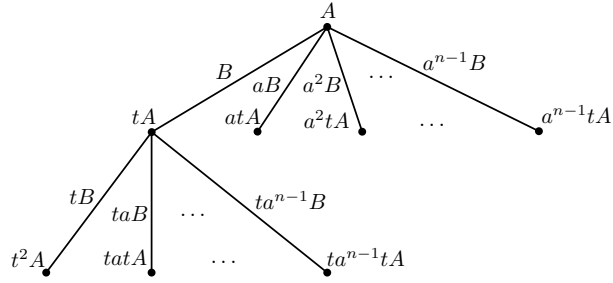}
\caption{Bass-Serre tree of $BS(1,n)$}\label{fig:BS1n-tree}
\end{figure}

The orientation on $T$ induced by the surjection $\varphi:G\rightarrow \mathbb{Z}$ allows us to distinguish two types of neighboors to any vertex $v$:

Let us look at the first of our two main examples: the Baumslag-Solitar group $BS(1,n)$. (The example of the lamplighter groups is postponed to the next section where we will also first prove that it can be seen as an HNN extension of type $(n,1)$). The Baumslag-Solitar group $BS(1,n)$ is an HNN extension for $H=A=\mathbb{Z}$, $B=n\mathbb{Z}$ and $\varphi:\mathbb{Z}\rightarrow n\mathbb{Z}$ given by multiplication by $n$,
$$BS(1,n)=\langle a,t\mid tat^{-1}=a^n\rangle.$$
Its Bass-Serre tree is depicted in Figure \ref{fig:BS1n-tree}.

\begin{lem}
\label{lemma: transitivity}
Let $G$ be an HNN extension such that $A=H$ and $B$ is a normal subgroup of $H$ of odd prime index $p=2k+1$.
Let $g \in G$ be an elliptic element. For any vertex $v$  of the Bass-Serre tree $T$ either $g(v)=v$ or the $p=2k+1$ vertices
$$ g^{-k}(v),\dots, g^{-1}(v),v,g(v),\dots, g^k(v)$$
are distinct.
\end{lem}

\begin{proof} Let $a\in A=H$ be any element not in the kernel of the natural surjection $A\rightarrow A/B\cong \mathbb{Z}_p$. Then $A=\sqcup_{j=-k}^k a^j B$. In the Bass-Serre tree of $G$, the $p$ direct descendants of the vertex $A$ are the vertices $a^{-k}tA,\dots, tA,\dots, a^k tA$ and are joined to $A$ through the edges $a^{-k}B,\dots, B,\dots, a^k B$ respectively. Observe that since $B$ is normal in $A$, any element $b\in B$ acts trivially on the direct descendants of the vertex $A$. Furthermore, $a$ and any of its powers $a^j$ where $p$ does not divide $j$ obviously acts cyclically on the first descendants of $A$.

By conjugation, we can suppose that our elliptic element is in fact $h=a^j b \in H=A$, with $b\in B$ and $-k\leq j \leq k$. If $j=0$ then $h$ acts trivially on the direct descendants of $A$, while if $j\neq 0$ then $h$ acts as a cyclic permutation of order $p$. This implies the lemma.
\end{proof}

The following Lemma is an immediate application of the classical ping-pong lemma for semigroups \cite[Proposition VII.2]{dlHBook} taking as ping-pong sets, the descendants of $x_iv$, for every $i$:

\begin{lem}[Ping-Pong Lemma]
Let $x_1,x_2,\ldots, x_r \in BS(1,p)$ act as positive hyperbolic automorphisms on the corresponding Bass-Serre tree $T$. Suppose that there exists a vertex $v\in T^0$ such that $\{x_1v,x_2v,\ldots,x_rv\}$ are leaves of a tree rooted at $v$. Then the set $\{ x_1,\ldots, x_r \}$ freely generates a free monoid.
\label{Key Lemma}
\end{lem}

\section{Lamplighter groups viewed as non-semidirect HNN extensions}\label{section: lamplighter groups}

The standard presentation for a restricted wreath product $G \wr \mathbb{Z}$ is also an HNN-extension, but the subgroups $A,B$ are both equal to $G$, so the corresponding Bass-Serre tree is a line, and the corresponding action of $G$ on a line is not useful for our goals. It has already been pointed out in \cite{Stalder-Valette} that there exists another HNN-extension presentation of any wreath product $G\wr \mathbb{Z}$ with indices $|G|$ and $1$ so that the corresponding Bass-Serre tree is a regular tree of valency $|G|+1$. For completeness, we include a proof of this fact for $\mathcal L_p=(\mathbb Z/ p\mathbb Z)\wr \mathbb{Z}$:

\begin{lem}
The lamplighter group $\mathcal L_p=(\mathbb Z/ p\mathbb Z)\wr \mathbb{Z}$ can be decomposed as an HNN- extension $H*_\theta$ with indices of the subgroups $[H:A]=1, [H:B]=p$.
\label{lem: non-semidirect}
\end{lem}

\begin{proof}
First, we find some useful presentation of the lamplighter group $\mathcal L_p$. Start with the standard presentation
\begin{equation*}
\mathcal{L}_p=\langle a, t \mid a^p=1, [t^m a t^{-m}, t^n a t^{-n} ]=1,\, m,n\in \mathbb{Z} \rangle,
\label{lamplighter-pres}
\end{equation*}
which can economically be rewritten as
\begin{equation}
\mathcal{L}_p=\langle a, t \mid a^p=1, [a, t^k a t^{-k} ]=1,\, k\in \mathbb{Z} \rangle
\label{lamplighter-pres2}
\end{equation}
since the element $[t^m a t^{-m}, t^n a t^{-n} ]$ can be obtained via conjugation of $[a, t^k a t^{-k} ]$ by a suitable power of $t$. Now we reduce the set of relations in \eqref{lamplighter-pres2} further to get
\begin{equation}
\mathcal{L}_p=\langle a, t \mid a^p=1, [a, t^k a t^{-k} ]=1,\, k\in \mathbb{N} \rangle,
\label{lamplighter-pres3}
\end{equation}
which is possible because $[a, t^{-k} a t^{k} ]$ is a consequence of a conjugate of $[a, t^k a t^{-k} ]$.

\smallskip

Now we show how to decompose the group $\mathcal L_p$ as a non-semidirect HNN-extension.
Consider the infinite direct sum  $D=\oplus_{\mathbb{N}_0} (\mathbb Z / p\mathbb Z)$ canonically generated by the set of elements $\{ a_0, a_1, a_2,\ldots \}$. The presentation of this group is
\begin{equation}
D=\langle a_0, a_1, a_2, \ldots \mid a_m^p=1,\, [a_m,a_n]=1,\, m,n\in \mathbb{N}_0 \rangle.
\end{equation}
Consider the HNN extension of $D$ given by the subgroups $H=D$ and $K=\langle a_1, a_2, \ldots \rangle$ and the isomorphism $f(a_i)=a_{i+1}$. Note that $[D:K]=p$. The HNN extension $D*_f$ then has the presentation
\begin{equation*}
D*_f=\langle t, a_0, a_1, a_2 \ldots \mid a_m^p=1,\,  [a_m,a_n]=1,\, ta_mt^{-1} = a_{m+1}, \; m,n\in \mathbb{N}_0  \rangle.
\end{equation*}
The relations $a_m^p=1$ with $m\geq 1$ can be excluded from this presentation because they follow from the relation $a_0^p=1$ and the series of relations $ta_mt^{-1} = a_{m+1}$. Then, repeatedly using the series $ta_mt^{-1} = a_{m+1}$ we substitute the letters $a_i$ in the commutators so that we get
$$
D*_f=\langle t, a_0, a_1, a_2 \ldots \mid a_0^p=1,\,  [t^m a_0 t^{-m},t^n a_0 t^{-n}]=1,\, t a_m t^{-1} = a_{m+1}, \; m,n\in \mathbb{N}_0  \rangle.
$$
Now we repeatedly remove the generators $a_m$ for all $m\geq 1$ and get the presentation
$$
D*_f=\langle t, a_0 \mid a_0^p=1,\,  [t^m a_0 t^{-m},t^n a_0 t^{-n}]=1,\; m,n\in \mathbb{N}_0  \rangle.
$$
Again, the relation $[t^m a_0 t^{-m},t^n a_0 t^{-n}]=1$ with $n\geq m$ follows from the relation $[a_0,  t^{n-m} a_0  t^{m-n}]=1$, so we obtain that
$$
D*_f=\langle t, a_0 \mid a_0^p=1,\,  [a_0,t^k a_0 t^{-k}]=1,\; k\in \mathbb{N}  \rangle,
$$
which is equivalent to the presentation \eqref{lamplighter-pres3} of the lamplighter group $\mathcal{L}_p$.
\end{proof}

It is quite obvious that the groups $\mathcal L_p$ tend to $\mathbb Z \wr \mathbb Z$ when $p$ tends to $\infty$. Actually, the following nice fact is also true:
\begin{prop}
The groups $BS(1,n)$ are factor groups of the wreath product $\mathbb Z \wr \mathbb Z$.
\label{L2-proposition}
\end{prop}

\begin{proof}
As seen above, the group $\mathcal L_{\mathbb Z}=\mathbb Z \wr \mathbb Z$ can be presented as
\begin{equation}
\mathbb Z \wr \mathbb Z=\langle a, t \mid [a, t^k a t^{-k} ]=1,\, k\in \mathbb{N} \rangle.
\label{lamplighterZ-pres}
\end{equation}

The presentations \eqref{Baumslag-Solitar-presentation} and \eqref{lamplighterZ-pres} prove the Proposition, since according to \eqref{Baumslag-Solitar-presentation}, for every positive $k$ the element $t^k a t^{-k}$ is a power of $a$, hence it commutes with $a$ so that the corresponding relation in \eqref{lamplighter-pres3} holds true.
\end{proof}

We will see later that $\lim_{p \to \infty}(\omega(BS(1,p)),\{a,t\})=1+\sqrt2=\omega(\mathbb Z \wr \mathbb Z,\{a,t\})$, which is some further evidence for the fact that $\mathbb Z \wr \mathbb Z$ is a limit of the groups $BS(1,n)$.

Now we can show that the classic lamplighter $\mathcal{L}_2$ gives the answer to Mann's question about growth of non-semidirect HNN-extensions (see \cite[Problem 1]{Mann}), proving a part of the Theorem 1, which we state as
\begin{prop} \label{part of THM 1} The minimal growth rate $\Omega(\mathcal L_2)$ of the lamplighter group $\mathcal L_2$ is realized on the generating set $\{a,t\}$ and it is equal to the golden ratio $\varphi=(1+\sqrt5)/2$.
\end{prop}
\begin{proof}
For the group $\mathcal{G} \wr {\mathbb{Z}}$ one can compute the exact growth series using the following formula of W.Parry from \cite[Corollary 3.3]{Parry}.
If $f_G(x)$ is the growth series of a finitely generated group $G$ then the growth series of $G \wr \mathbb{Z}$ can be obtained as
\begin{equation}
f_{G \wr \mathbb{Z}}(x)=\frac{f_G(x)(1-x^2)^2(1+x f_G(x))}{(1-x^2 f_G(x))^2(1-x f_G(x))}.
\label{Parry-formula}
\end{equation}

We use this formula to compute the growth series for $\mathcal{L}_2$.

$$
f_{\mathcal L_2}(x)=\frac{(1+x)(1-x^2)^2(1+x(1+x))}{(1-x^2(1+x))(1-x(1+x))}=
\frac{(1+x)(1-x^2)^2(1+x+x^2)}{(1-x^2+x^3)(1-x-x^2)}.
$$

The factors in nominator have roots on the unit circle, whilst the factors of the denominators give two roots inside the unit circle, whose reciprocals are the golden ratio $\varphi=(1+\sqrt5)/2$ and the so-called ``plastic number''\footnote{Notably $\psi=\Omega(GL(2,\mathbb{Z}))=\Omega(PGL(2,\mathbb{Z}))$, see \cite{Bucher-Talambutsa-IJM} for more information about this number.}. Since $\varphi>\psi$, we get $\omega(\mathcal{L}_2,\{a,t\})=\varphi$. As $\mathcal{L}_2$ is a non-semidirect HNN extension due to Lemma~\ref{lem: non-semidirect}, we may apply the Theorem 1 from \cite{Mann} to get the lower bound $\omega(\mathcal{L}_2)\geqslant\varphi$ and finally conclude that $\omega(\mathcal{L}_2)=\varphi$.
\end{proof}

The equality $\omega(\mathcal L_2,\{a,t\})=\varphi$ was also mentioned in \cite[p.1997]{LPP} by Lyons-Pemantle-Peres, and follows from the observation that there is a subtree in the Cayley graph of $\mathcal{L}_2$ which is a Fibonacci tree.

\begin{rem}
It would be interesting to find a natural (maybe geometric) reason for the group $\mathcal{L}_2$ to have the ``second biggest growth rate'' equal to the plastic number $\psi$.
\end{rem}

\section{Growth rates computations and estimates}

We collect in this section some explicit computations and estimates on growth rates. Lemma \ref{lem: lengthMonoid}, which is proved in \cite[Lemma 6]{Bucher-Talambutsa-IJM}, will be used extensively in the proofs of Theorems  \ref{Main Theorem} and \ref{Main Theorem 2} in combination with our Ping-Pong Lemma \ref{Key Lemma}. The exact growth rates of some Baumslag-Solitar groups and lamplighters groups are computed in Lemma \ref{LemmaGrowthRatesOfBS} and the last Lemma \ref{lem:RootsComparison} allows us to compare some particular roots.

\begin{lem} \label{lem: lengthMonoid} Let $G$ be a group generated by a
finite set $S$. Suppose that there exists a set $\{x_1,\dots,x_k\}\subset G$ generating a
free monoid inside $G$. Set $\ell_i=\ell_{G,S}(x_i)$, for $i=1,\ldots,k$,
and $m=\max\{ \ell_1,\dots,\ell_k\}$. Then $\omega(G,S)$ is greater or equal
to the unique positive root of the polynomial
\begin{equation}
 Q(z)=z^m-\sum_{i=1}^k z^{m-\ell_i}.
 \label{char-polynomial}
 \end{equation}
\label{lem:Computation for monoid}
\end{lem}

As mentioned in the introduction we can easily compute the growth rate of the lamplighters and Baumslag-Solitar group with respect to the canonical generators from the growth functions found by Parry \cite{Parry} and Collins, Edjvet and Gill \cite{Collins-Edjvet-Gill} respectively. Recall that for any integer $k \geq 1$ we consider the polynomial
\[T_k(x)=x^{k+1}-x^k-2x^{k-1}-\ldots-2x-2.\]
Due to Descartes rule of signs, $T_k$ has single positive root, which we denote by $\omega_k$.

\begin{lem}
\noindent (a) The growth rate $\omega(\mathcal{L}_2,\{a,t\})$ is equal to $\frac{1+\sqrt5}2$.

(b) For any $k \geq 1$ we have that
\[
\omega(BS(1,2k+1),\{a,t\})=\omega(\mathcal{L}_{2k+1},\{a,t\})=\omega_k,
\]

\noindent (c) The growth rate $\omega(BS(1,2),\{a,t\})$ is equal to the positive root of $x^3-x^2-2$.
\label{LemmaGrowthRatesOfBS}
\end{lem}

\begin{proof} (a) This part was already shown in the Proposition \ref{part of THM 1}.

(b) Another elegant formula by Parry (see \cite[Theorem 4.1]{Parry}) allows to compute the growth rate of the wreath product $G\wr \mathbb{Z}$. If $S$ is a finite generating set for the group $G$ then $\omega(G\wr \mathbb{Z},S\cup \{t\})=1/\kappa$, where $\kappa$ is the smallest positive zero of the function $1-xf_{G,S}(x)$.
Taking $f_{\mathbb{Z}/(2k+1)\mathbb{Z},\{a\}}(x)=1+2x+2x^2+\ldots+2x^{k-1}$ we get that $\omega(\mathcal{L}_{2k+1},\{a,t\})=1/\kappa_k$, where $\kappa_k$ is the smallest positive zero of the polynomial $R_k(x)=1-x-2x^2-\ldots-2x^{k+1}$. The polynomials $R_k$ and $T_k$ are reciprocal, so indeed we get that $\omega(\mathcal{L}_{2k+1},\{a,t\})=1/\omega_k$.

To prove that $\omega(BS(1,2k+1),\{a,t\})=\omega_k$ we use the following explicit formula from \cite{Collins-Edjvet-Gill}, which gives a power series $\Sigma_k(x)=\sum_{m=0}^{\infty} f(m) x^m$ for the growth function $f(m)=f_{BS(1,n),\{a,t\}}(m)$. For the case $n=2k+1$ they obtain
\begin{equation}
\Sigma_n(z)=\frac{(1+x^2-2x^{k+2})(1+x-2x^{k+2})(1+x)^2(1-x)^3}{(1-x-x^2-x^3+2x^{k+3})^2(1-2x-x^2+2x^{k+2})}.
\label{Collins-Edjvet-Gill-formula}
\end{equation}
Then the growth rate $\omega(BS(1,2k+1),\{a,t\})$ is equal to $1/\alpha$, where $\alpha$ is the smallest positive pole of the function $\Sigma_{n}(x)$. Since $1<\omega(BS(1,2k+1),\{a,t\})<3$, we have bounds $1/3<\alpha<1$.
We will first prove that $\alpha=\gamma_2$, where $\gamma_2$ is the smallest positive root of the second factor $Q_2(x)=1-2x-x^2+2x^{k+2}$ of the denominator of \eqref{Collins-Edjvet-Gill-formula}. Let $\gamma_1$ be the smallest positive root of the first factor $Q_1(x)=1-x-x^2-x^3+2x^{k+3}$. Note that $Q_1(0)=Q_2(0)=1$ and $Q_1(1)=Q_2(1)=0$, so the numbers $\gamma_1,\gamma_2$ are well defined and $0<\gamma_1,\gamma_2\leq 1$.

Since the difference function
\[
Q_1(x)-Q_2(x)=x-x^3+2x^{k+2}-2x^{k+3}=x(1-x^2)+2x^{k+1}(1-x)
\] is non-negative on $[0,1]$, we obtain that $\gamma_1\geq \gamma_2$.

To show that $\alpha=\gamma_2$ we are left to prove that $\gamma_2$ is not a root of the nominator. The factors $(1+x)^2$ and $(1-x)^3$ do not have roots on the interval $I=(1/3,1)$, and we will check that $P_1(x)=1+x^2-2x^{k+2}$ and $P_2(x)=1+x-2x^{k+2}$ have no common roots with $Q_2(x)$ on $I$. This is true, since otherwise either $Q_2(x)+P_1(x)=2-2x$ or $Q_2(x)+P_2(x)=(2+x)(1-x)$ would have a root on $(1/3,1)$, which is false.

We can factorize $Q_2(x)$ as $(1-x)Z(x)$ with $Z(x)=1-x-2x^2-\ldots-2x^{k+1}$. Since the polynomial $Z(x)$ is reciprocal to the polynomial $T(x)$ from the statement, the part (b) of Lemma is proved.

(c) Here we use another formula from \cite{Collins-Edjvet-Gill} that is
\[
\Sigma_2(x)=\frac{(1-x)^2(1+x)^2H(x)}{(1-x-2x^3)(1-x^2-2x^5)^2},
\]
where $H(x)=1+3x+8x^2+12x^3+16x^4+20x^5+22x^6+16x^7+14x^8+12x^9+4x^{10}$.

We follow the same strategy as in the part (b), and first make sure that the positive root of the polynomial $Q_1(x)=1-x-2x^3$ is smaller than the one of $Q_2(x)=1-x^2-2x^5$, because $Q_2(x)-Q_1(x)=x(1-x)+2x^3(1-x^2)>0$ on $(0,1)$. Then, making tedious computations or using a computer, one gets that $\mathrm{GCD}(H(x),Q_1(x))=1$, so the smallest pole of $\Sigma_2(x)$ indeed comes from $Q_1(x)$. Again, $Q_1(x)$ is reciprocal to $x^3-x^2-2$, and the part (c) is also proved.
\end{proof}

The next lemma will allow us to compare $\omega_k$ with the growth rate of some free monoid in the proof of Theorem \ref{Main Theorem 2}.

\begin{lem}
Let $k\geq 1$ be an integer and $\delta_k$ be the unique positive root of the polynomial $D_k(x)=x^{2k+1}-2x^{2k}-2x^{2k-2}-\ldots-2x^2-2$. Then
\[\frac{1+\sqrt5}{2}\leq \omega_k \leq \delta_k<1+\sqrt2.\]
\label{lem:RootsComparison}
\end{lem}
\begin{proof}
The inequality $(1+\sqrt5)/2\leq \omega_k$ may be proven directly, but actually we already know that $\omega(BS(1,2k+1),\{a,t\})=\omega_k$ and $\Omega(BS(1,2k+1))\geq (1+\sqrt5)/2$ as proved by Mann.

Since $T_k(1),P_k(1)<0$ and $T_k(+\infty)=P_k(+\infty)=+\infty$ we get $\delta_k, \omega_k>1$. Consider the polynomials $D(x)=(x^2-1)D_k$ and $T(x)=(x^2-1)T_k=(x+1)(x-1)T(x)$. After a simple calculation we get
\begin{eqnarray*}
D(x)&=&
x^{2k+3}-2x^{2k+2}-x^{2k+1}+2,\\
T(x)&=&
x^{k+3}-x^{k+2}-3x^{k+1}-x^k+2x+2.
\end{eqnarray*}
As $(x^2-1)>0$ on $(1,+\infty)$ and $D(1+\sqrt2)=2>0$, we get that $\delta_k\in(1,1+\sqrt2)$.

Since $T(1)=D(1)=0$ and $T(1+\varepsilon),D(1+\varepsilon)>0$ for small $\varepsilon$, in order to show the inequality $\omega_k \leq \delta_k$ it suffices to show that $T(x)\geq D(x)$ on the interval $(1,1+\sqrt2)$.

Consider the difference function
\begin{equation*}
\begin{aligned}
D(x)-T(x)&=x^{2k+3}-2x^{2k+2}-x^{2k+1}-x^{k+3}+x^{k+2}+3x^{k+1}+x^k-2x\\
&=(x^k-1)(x^{k+1}-1)(x^2-2x-1)-(x^2-1).
\end{aligned}
\end{equation*}
Since the polynomials $x^k-1$ and $x^{k+1}-1$ are positive on $(1,+\infty)$ and $x^2-2x-1$ is negative on $(1,1+\sqrt2)$,
we indeed have that $D(x)-T(x)<0$ on $(1,1+\sqrt2)$, which proves the lemma.
\end{proof}

\section{Proofs of Theorems \ref{Main Theorem} and \ref{Main Theorem 2}}

\begin{proof}[Proof of theorem \ref{Main Theorem 2}]  Let $G=H*_\theta$ be an HNN extension relative to an isomorphism $\theta:A\rightarrow B$ with $A=H$ and $B$ a normal subgroup of prime index $p$ in $H$. Let $S$ be any generating set for $G$. We need to show that $\omega(G,S)\geq (1+\sqrt{5})/2$ for $p=2$ and $\omega(G,S)\geq \omega_k$ for $p=2k+1$.

As explained above (see Section \ref{section HNN}), the natural surjection $\varphi:G\rightarrow \mathbb{Z}$ ensures the existence of a hyperbolic element in $S$. Furthermore, upon replacing $x$ by $x^{-1}$ we can suppose that $x$ is a positive element. Since the action of $G$ is transitive on its $(p+1)$-regular Bass-Serre tree, there exists an element in $S$ not preserving the axis $L_x$ of $x$. We distinguish two cases according to this element being elliptic or hyperbolic.

\smallskip

\noindent\textsc{Case 1 (elliptic)} There exists an elliptic element $z\in S$ such that $z(L_x)\neq L_x$.

For $p=2$, we consider the set
\[
M=\{x, \ zx \},
\]
while for odd primes $p=2k+1$,
\[
M=\{x,\,zx,z^2x, \ldots , z^kx,\,\,z^{-1}x,z^{-2}x, \ldots , z^{-k}x\}.
\]
In either cases, we will show that $M$ freely generates a free monoid.

Since any vertex has only one direct ascendant, if a vertex is in the fixed point set of $z$, then all its ascendants are. For the same reason, any two ascending rays meet, so there exists a vertex of the axis of $x$ which is fixed by $z$. Let $v$ be the lowest vertex on $L_x \cap \mathrm{Fix}(z)$. Then $x(v)$ is a descendant of $v$, which is not in the set $\mathrm{Fix}(z)$, hence the vertices
\[x(v),\, zx(v), \quad \mathrm{for\ } p=2,\]
and by Lemma \ref{lemma: transitivity}, the vertices
\[x(v),\,\,zx(v), \ldots , z^kx(v),\,\,z^{-1}x(v), \ldots , z^{-k}x(v), \quad \mathrm{for \ odd\ } p=2k+1,\]
are all distinct leaves of a tree rooted at $v$, so $M$ freely generates a free monoid due to the Ping-Pong Lemma \ref{Key Lemma}. Lemma \ref{lem:Computation for monoid} now implies that $\omega(G,S)$ is greater or equal to the unique positive root of
$$ z^2-2z-1, \quad \mathrm{for\ } p=2,$$
which is precisely the golden ratio $(1+\sqrt{5})/2$, while for $p=2k+1$, it is greater or equal to the unique positive root of
$$ T_k(z)=z^{k+1}-z^k-2z^{k-1}-\dots -2z-2, $$
which is $\omega_k$ by definition.

\smallskip

\noindent\textsc{Case 2 (hyperbolic).} There exists a hyperbolic element $y\in S$ such that $y(L_x)\neq L_x$. Upon replacing $y$ by its inverse, we can suppose that $y$ is a positive hyperbolic. Since $y$ preserves its axis $L_y$, this implies that the axes $L_x$ and $L_y$ are different. This already implies that $\omega(BS(1,p),S)\geq 2$ (see \cite[Lemma]{dlHBucher} or Lemma \ref{lem:Computation for monoid} with $\ell_1=\ell_2=1$). Since for $p=2,3$ we have
\[\omega(BS(1,2),\{a,t\})<\omega(BS(1,3),\{a,t\})=2,\]
we can suppose that $p\geq 5$, and again $p=2k+1$.

We consider four subcases, according to the situations when \\
{\textsc a.} $\ell(x)=\ell(y)$, {\textsc b.} $2\ell(y)<\ell(x)$, {\textsc c.} $\ell(x)=2\ell(y)$ and {\textsc d.} $\ell(y)<\ell(x)< 2\ell(y)$.

\smallskip

\noindent\textsc{Case 2a.} $\ell(x)=\ell(y)$. Note that the element $yx^{-1}$ is elliptic and $yx^{-1}(L_x)\ne L_x$. We can apply the claim of \textsc{Case 1} to $x$ and $z=yx^{-1}$ to conclude that the set
\[
\{x,y,yx^{-1}y,\ldots, (yx^{-1})^{k-1}y,xy^{-1}x,\dots, (xy^{-1})^kx\}
\]
freely generates a free monoid. Then Lemma \ref{lem:Computation for monoid} shows that $\omega(BS(1,2k+1),S)\geq \delta_k$,
where $\delta_k$ is the single positive root of the polynomial $D_k(x)=x^{2k+1}-2\sum_{m=0}^{k}x^{2m}$. Finally, Lemma \ref{lem:RootsComparison} gives the desired inequality $\omega(BS(1,2k+1)) \geq \delta_k \geq \omega_k$.

We can now suppose that $\ell(y)<\ell(x)$ and distinguish three further subcases:

\textsc{Case 2b.} $ 2\ell(y)<\ell(x)$

We will show that the infinite family
$$\{ y^{-2}x, y^{-1}x, x, yx,y^2x,\dots, y^sx,\dots, \\
   yx^{-1}yx,y^2x^{-1}yx,\dots, y^sx^{-1}xy,\dots \}$$
which is maybe better described as
$$\{ y^sx\mid s\geq -2\} \cup \{ y^sx^{-1}yx\mid s\geq 1\}$$
freely generates a free monoid. Then, taking as free generators only the $2k+1$ elements
\[x, yx,y^2x,\dots, y^kx, y^{-1}x,y^{-2}x,
yx^{-1}yx,y^2x^{-1}yx,\dots, y^{k-2}x^{-1}xy\] we get that $\omega(G,S)$ is by Lemma \ref{lem:Computation for monoid} greater or equal to the unique positive root of
$$ T_k(z)=z^{k+1}-z^k-2z^{k-1}-\dots -2z-2, $$
which is $\omega_k$ by definition.

To prove that the above infinite family freely generates a monoid, let $v_0$ be the lowest vertex on $L_x\cap L_y$ and let $v_x\in L_x$ and $v_y\in L_y$ be the corresponding direct descendants of $v_0$. We aim at applying the Ping-Pong Lemma \ref{Key Lemma} to the vertex $w=x^{-1}(v_x)$, see Figure \ref{fig:Case2Y<X}.

First notice that since $v_x\notin L_y$, the translates $y^s x (w)=y^s (v_x)$ are all distinct, branching from $L_y$ at $y^s(v_0)$. Furthermore, for $-2\leq s$, the highest such translate is $y^{-2}x(w)=y^{-2}(v_x)$ which is strictly below $y^{-2}(v_0)$ by construction. Now $w=x^{-1}(v_x)$ is equal or above $y^{-2}(v_0)$ since $2\ell(y)<\ell(x)$. This already implies that the infinite subfamily $\{ y^sx\mid -2\leq s \}$ freely generates a free monoid.

Second consider the vertex $y(v_x)$. It is branching from $L_x$ at $v$ and the first vertex from $L_x\cap L_y$ to $y(v_x)$ is $v_y$. It follows that $x^{-1}y(v_x)$ does not belong to $L_x$ either and is branching at $x^{-1}(v)$ from $L_x$ and hence also from $L_y$. It follows that all the translates $y^sx^{-1}yx(w)=y^sx^{-1}y(v_x)$ belong to different branches of $L_y$, branching at $y^sx^{-1}(v_0)$. Since $\ell(y)\geq 1$, for $1\leq s$ the branch points are below or equal to $w=x^{-1}(v_x)$.

\begin{figure}
\includegraphics[width=125mm]{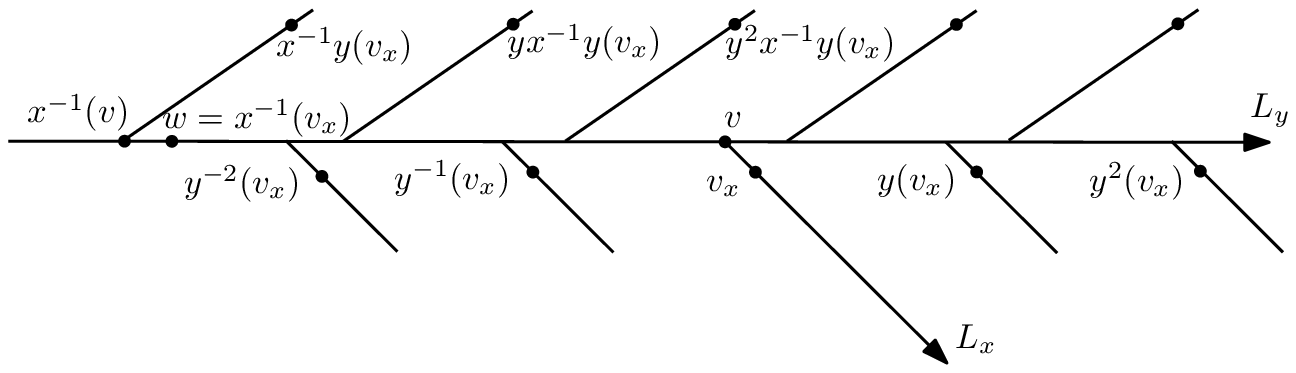}
\caption{Case $ 2\ell(y)<\ell(x)$.}
\label{fig:Case2Y<X}
\end{figure}

If $\ell(x)$ is not a multiple of $\ell(y)$ the two families of branching points are different and we are done. If $\ell(x)=m\ell(y)$ for some $m>2$ we need to check that $y^{n+m} x^{-1}(v_y)\neq y^n v_x$ and it is enough to check it for $n=0$. Consider the elliptic element $y^mx^{-1}$. It fixes $v_0$, sends $v_x$ to $v_y$ and $v_y$ to $y^mx^{-1}(v_x)$ which cannot be equal to $v_x$ otherwise the action on the direct descendants of $v_0$ of the elliptic element $y^mx^{-1}$ would not be transitive, contradicting Lemma \ref{lemma: transitivity}.

\textsc{Case 2c.} $\ell(x)=2\ell(y)$.

It is enough to show that the set
$$\{ x,y,xy^{-1}x,xy^{-2}x,xy^{-1}xy^{-1}x,y^2x^{-1}y,xyx^{-1}y\}$$
freely generates a free monoid. Then, using Lemma \ref{lem:Computation for monoid} we get that $\omega(BS(1,k))$ is at least $\gamma$, where $\gamma$ is the root of the polynomial $F(x)=x^5-2x^4-x^2-3x-1$. Since $F(x)=(x^2-2x-1)(x^3+x+1)$, we get that $\gamma=1+\sqrt2$, and again Lemma \ref{lem:RootsComparison} gives
the desired inequality $\omega(G,S)\geq \omega_k$.

Let as above $v$ be the lowest vertex on  $L_x\cap L_y$. We aim at applying the Ping-Pong Lemma \ref{Key Lemma} to the vertex $v$. Let $v_x\in L_x$ and $v_y\in L_y$ be the corresponding direct descendants of $v_0$.

The elliptic transformation $b=y^2x^{-1}$ fixes $v$ and takes $v_x$ to $v_y$. Thus its action on the direct descendants of $v$ is nontrivial and hence transitive. Since we assume $p\geq 4$, it follows by Lemma \ref{lemma: transitivity} that the image $v_+=y^2x^{-1}(v_y)$ of $v_y$ and the preimage $v_-:=xy^{-2}(v_x)$ of $v_x$ give four distinct direct descendants of $v_0$ as depicted in Figure \ref{fig:Case2Y=X}.

\begin{figure}
\includegraphics[width=70mm]{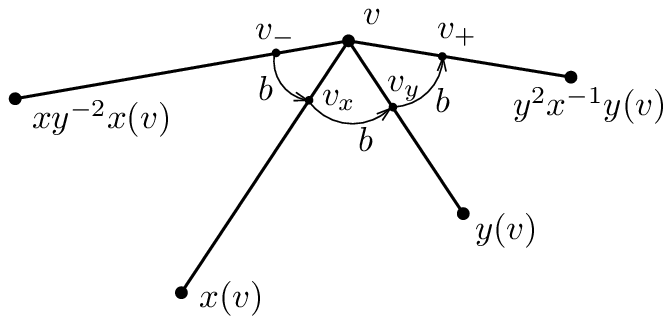}
\caption{Case $\ell(x)=2\ell(y)$: the action of the elliptic element $b=y^2x^{-1}$.}
\label{fig:Case2Y=X}
\end{figure}

Observe that $y^2x^{-1}y(v)$ is on the branch through $v$ and $v_+$, while $xy^{-2}x(v)$ is on the branch through $v_0$ and $v_-$. Thus the four elements $xv,yv,xy^{-2}x(v)$ and $y^2x^{-1}y(v)$ have distinct geodesics to $v$.

We now forget about $xy^{-2}x(v)$ and look at the image of the tree rooted at $v$ of the three remaining elements through the hyperbolic transformation $xy^{-1}$. The root $v$ is mapped on the segment from $v$ to $x(v)$. The vertex $y(v)$ is mapped to $x(v)$, and the two remaining leaves are sent to vertices branching from $L_x$ at $xy^{-1}(v)$.

Iterating this procedure but only on $xy^{-1}(v), x(v)$ and $xy^{-1}x(v)$ shows that $xy^{-1}xy^{-1}x(v)$ is branching from the segment between $xy^{-1}(v)$ and $xy^{-1}x(v)$. We have thus proven that the seven vertices are leaves of a tree rooted at $v$, as illustrated in Figure \ref{fig:DetailedCase2Y=X}, which finishes the proof of this case.

\begin{figure}
\includegraphics[width=100mm]{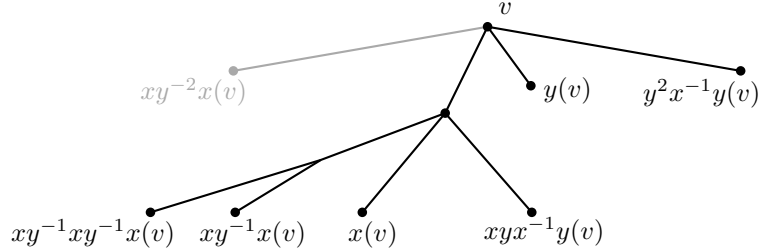}
\caption{Case $\ell(x)=2\ell(y)$: The subtree to which we apply the Ping-Pong Lemma \ref{Key Lemma}.}
\label{fig:DetailedCase2Y=X}
\end{figure}

\textsc{Case 2d.} $\ell(y)<\ell(x)<2\ell(y)$.

We will show that the set
$$\{ x,y,xy^{-1}x,xy^{-2}x,yx^{-1}y\}$$
freely generates a free monoid. Since the corresponding polynomial $x^4-2x^3-2x-1=x(x^2+1)(x^2-2x-1)$ has only one positive root $1+\sqrt 2$, this will prove this case.

Set $a=\ell(x)$ and $b=\ell(y)$. The proof decomposes in the two cases $b<a\leq (3/2)b$ and $(3/2)b\leq a<2b$ with an additional small argument needed in the equality case.

In case $b<a\leq (3/2)b$ we aim at applying the Ping-Pong Lemma \ref{Key Lemma} to the vertex $w=xy^{-2}(v)$. (See Figure \ref{fig:CaseY<X<3/2Y}.) This vertex is on the intersection of the axes $L_x\cap L_y$ at distance $2b-a$ above $v$. Of the five images of $w$, only $x(w)$ is on the axis $L_x$, at  distance $a$ below $w$ and hence $2(a-b)$ below $v$. The four other images are not in $L_x$ and we will determine their projection on $L_x$.

The image $y(w)$ is on the axis $L_y$ at distance $b$ below $w$ and hence at distance $a-b$ from its projection $v\in L_x$. Since the axis of the hyperbolic transformation $xy^{-2}$ contains $L_x\cap L_y$ and at least the vertex $v_y\in L_y$, the segment $[v,x(w)]$, which intersects $L_{xy^{-2}}$ only at $v$ is mapped by $xy^{-2}$ to the segment $[w,xy^{-2}x(w)]$ which intersect $L_{xy^{-2}}$ and hence $L_x$ only in $w$. Similarly, the axis of $xy^{-1}$ contains $L_x\cap L_y$ and at least the vertex $v_x\in L_x$, so that the hyperbolic transformation $xy^{-1}$ takes the segment $[v,x(v)]$ to the segment $[xy^{-1}(v),xy^{-1}x(v)]$ which intersects $L_{xy^{-1}}$ and hence $L_x$ precisely in $xy^{-1}(v)$ which is at distance $a-b$ from both $v$ and $x(v)$. Finally, the axis of $yx^{-1}$ contains $L_x\cap L_y$ and at least the vertex $v_y\in L_y$, so that applying $yx^{-1}$ to the segment $[v,y(w)]$ we obtain the segment $[yx^{-1}(v),yx^{-1}y(w)]$ which intersects $L_x$ in $yx^{-1}(v)$ which is at distance $a-b$ above $v$ and hence at distance $3b-2a\geq 0$ below $w$. If the inequality is strict, the claim immediately follows from the Ping-Pong Lemma \ref{Key Lemma}. If $3b-2a=0$, we will see below how to show that the segments $[yx^{-1}(v),yx^{-1}y(w)]$ and $[w,xy^{-2}x(w)]$ only intersect at $w=yx^{-1}(v)$.


\begin{figure}
\includegraphics[width=75mm]{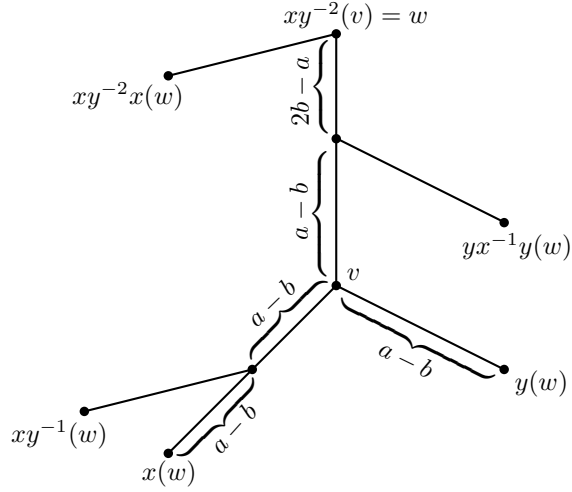}
\caption{Case $b<a<3/2b$.}
\label{fig:CaseY<X<3/2Y}
\end{figure}

If $(3/2)b\leq a<2b$ the argument is completely analogous, except that the vertex $yx^{-1}(v)$ is above or equal to $w=xy^{-2}(v)$. Thus we want to replace $w$ by $w':=yx^{-1}(v)$ and apply the Ping-Pong Lemma \ref{Key Lemma} to this vertex $w'$. (See Figure \ref{fig:Case3/2Y<X<2Y}. This vertex is on the intersection of the axes $L_x\cap L_y$ at distance $a-b$ above $v$. Of the five images of $w'$, only $x(w')$ is on the axis $L_x$, at  distance $a$ below $w$ and hence $b$ below $v$. The four other images are not in $L_x$ and we will determine their projection on $L_x$.

The image $y(w')$ is on the axis $L_y$ at distance $b$ below $w$ and hence at distance $2b-a$ from its projection $v\in L_x$. For the three other image points, the proof is identical to the above case, replacing $w$ by $w'$.

\begin{figure}
\includegraphics[width=75mm]{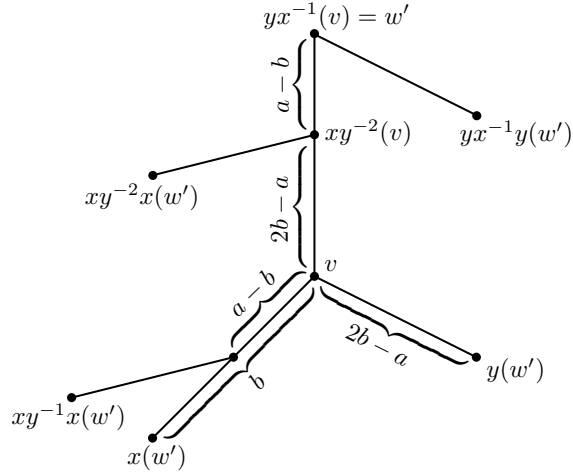}
\caption{Case $3/2b<a<2b$.}
\label{fig:Case3/2Y<X<2Y}
\end{figure}

In the equality case the two vertices $w=w'$ agree. Let $v_1$, respectively $v_2$ be the first vertex after $w$ on the geodesic to $xy^{-2}(w)$, respectively $yx^{-1}y(w)$. We need to show that $v_1\neq v_2$. Let $v_a$ be the direct descendant of $w$ on the geodesic to $v$. The ordered pair $(v_1,v_a)$ is mapped to $(v_x,v_y)$ by $y^2x^{-1}$, which are further mapped to $(v_a,v_2)$ by $yx^{-1}$. Thus the elliptic element $yx^{-1}y^2x^{-1}$ sends the ordered pair $(v_1,v_a)$ to $(v_a,v_2)$ and since $p\geq 3$ and elliptic elements act either trivially or transitively on direct descendants of a fixed point by Lemma \ref{lemma: transitivity} it follows that $v_1\neq v_2$, which finishes the proof of this case and of the theorem. \end{proof}



\begin{proof}[Proof of theorem \ref{Main Theorem}]  In view of Lemma \ref{LemmaGrowthRatesOfBS}, Theorem \ref{Main Theorem} follows immediately from Theorem \ref{Main Theorem 2} except in the case of $BS(1,2)$ where we need a better understanding of its action on the Bass-Serre to obtain the accurate lower bound of $\omega(BS(1,2),\{a,t\})=\beta$, where $\beta$ is the unique real root of $x^3-x^2-2$.

Let $S$ be a generating set for $BS(1,2)$. As in the proof of Theorem \ref{Main Theorem 2}, the case where $S$ contains two hyperbolic elements with different axis immediately gives the lower bound of $\omega(BS(1,2),S)\geq 2>\beta$. We thus only have to treat the corresponding elliptic case, that is, there exists a positive hyperbolic element $x\in S$ with axis $L_x$ and an elliptic element $z\in S$ such that $z(L_x)\neq L_x$.

As observed in the elliptic case of the proof of Theorem \ref{Main Theorem 2} the intersection of $L_x$ with the fixed point set of $z$ is nonempty. Upon conjugating the generating set $S$, we can suppose that the lowest vertex on $L_x$ fixed by $z$ is $A$, which implies that $z$ belongs to $A$. Since $z$ does not fix the direct descendants $tA$ and $atA$ it must be an odd power of $A$.

Consider the action of $a$ on the second generation of descendants of $A$, that is $t^2A, tatA, at^2A$ and $atatA$. The action has order four, mapping $t^2A\mapsto at^2A\mapsto a^2t^2A=tatA\mapsto atatA\mapsto a^2tatA=t^2A$. The action of $z$, as an odd power of $A$ is thus necessarily equal to the action of $a$ or $a^{-1}$ on these second generation descendants. It follows that $xA, zx^2A$ and $z^{-1}x^2A$ are leaves of a tree rooted at $A$, and hence $x,zx^2,z^{-1}x^2$ generate a free monoid by the Ping-Pong Lemma \ref{Key Lemma}.  Since these elements have lengths $1,3$ and $3$ respectively, we can invoke \ref{lem:Computation for monoid} to conclude that the grow rate of $BS(1,2)$ with respect to $S$ is greater or equal to the greatest and unique real root of $x^3-x^2-2$. Finally, Lemma \ref{LemmaGrowthRatesOfBS} gives
\[
\omega(BS(1,2),S)\geq \omega(BS(1,2),\{a,t\}),
\]
which finishes the proof of the theorem.
\end{proof}

The next lemma will be needed to prove Corollary \ref{cor: ZwrZ}.

\begin{lem}
The limit $\lim\limits_{k \to \infty} \omega_k = 1+\sqrt 2$ exists.
\label{lem: limitofomegas}
\end{lem}
\begin{proof}
From Lemma~\ref{lem:RootsComparison} and the definition of $\omega_k$ we know that $\omega_k$ is a single positive root of the polynomial $T_k(x)$, and $(1+\sqrt5)/2<\omega_k<1+\sqrt 2$ for every $k\geq 1$. Then the reciprocal polynomial $R_k(x)=1-x-2x^2-\ldots-2x^k-2x^{k+1}$ has a single positive root $1/\omega_k$ which belongs to the interval $I=(1/3,2/3)$. Consequently the polynomial
\[
R'_k(x)=(1-x)R_k=(1-x)^2-2x^2(1-x^k)=1-2x-x^2+2x^{k+2}
\]
also has two positive roots: $1$ and $1/\omega_k$. Obviously, for $k\to\infty$ the polynomials $2x^{k+2}$ uniformly converge to the zero function on the enlarged interval $I'=(1/4,3/4)$. For this reason the roots $1/\omega_k$ of $R'_k(x)$ on $I$ converge to the root of the polynomial $1-2x-x^2$ on $I$, and the latter root is equal to $\sqrt2-1=1/(1+\sqrt2)$, which proves the lemma.
\end{proof}

\begin{proof}[Proof of Corollary \ref{cor: ZwrZ}]

We use Parry's formula \eqref{Parry-formula} to compute the series $\Sigma(x)$ for the growth function $\mathbb Z \wr \mathbb Z$ with respect to the generating set $\{a,t\}$:
$$
\Sigma_{}(x)=
\frac{(1-x^2)^3(1+x^2)}{(1-x-x^2-x^3)^2(1-2x-x^2)}=\frac{(1+x)^2(1-x^2)^3(1+x^2)}{(1-x^4)^2(1-2x-x^2)}.
$$
All the roots of the nominator and the denominator lie on the unit circle except for the roots of $1-2x-x^2$. The reciprocal of the smallest root is equal to $\sqrt2+1$, hence this is the value for $\omega(\mathbb Z \wr \mathbb Z,\{a,t\})$.

\smallskip

Now we will show that $\Omega(\mathbb Z \wr \mathbb Z)=1+\sqrt2$. We already know that $\Omega(\mathbb Z \wr \mathbb Z )\leq 1+\sqrt2$. Suppose that $\omega(\mathbb Z \wr \mathbb Z)= 1+\sqrt2-\varepsilon$, where $\varepsilon>0$. As any group $\mathcal L_p$ is a factor group of the group $\mathbb Z \wr \mathbb Z$, then for any prime $p$ we have $\omega(\mathcal{L}_p)\leq 1+\sqrt2-\varepsilon$ which contradicts the limit equality $\lim\limits_{p\to \infty}\omega(\mathcal{L}_p)=\lim\limits_{k \to \infty} \omega_k=1+\sqrt 2$ proven in Lemma~\ref{lem: limitofomegas}.
\end{proof}

\bibliographystyle{amsalpha}

\end{document}